\documentclass{article}

\usepackage[T2C]{fontenc}
\usepackage[cp1251]{inputenc}
\usepackage[english]{babel}

\usepackage[tbtags]{amsmath}
\usepackage{amsfonts,amssymb,mathrsfs,amsmath,amscd}

\newtheorem{theorem}{Theorem}

\newtheorem{llemma}{Lemma}

\newtheorem{definition}{Definition}

\begin{document}

\title{The Discontinuity Group\\of a Locally Bounded Homomorphism\\%
of a Lie Group into a Lie Group\\Is Commutative}

\author{\leftline{A.~I.~Shtern$^{*,**,***,1}$}\\
\leftline{\hskip-60pt$^*$Moscow Center for Fundamental and Applied Mathematics,
Moscow, 119991
Russia}\\
\hskip-60pt$^{**}$Department of Mechanics and Mathematics, Moscow State
University,
Moscow, 119991 Russia\\
\hskip-60pt$^{***}$Scientific Research Institute for System Analysis of
the Russian Academy of Sciences\\\hskip-220pt(FGU FNTs NIISI RAN), Moscow, 117312 Russia\\
\hskip-60pt E-mail: {\tt $^1$aishtern@mtu-net.ru}}

\date{July 20, 2023}

\maketitle
\begin{abstract}
We prove that the discontinuity group of every locally bounded homomorphism of
a Lie group into a Lie group is not only compact and connected, which is known,
but is also commutative.
\end{abstract}

\section{Introduction}

Recall that a (not necessarily continuous) homomorphism $\pi$ of a topological
group $G$ into a topological group~$H$ is said to be \emph{relatively
compact\/} if there is a neighborhood $U=U_{e_G}$ of the identity element $e_G$
in~$G$ whose image $\pi(U)$ has compact closure in~$H$.

Obviously, a homomorphism into a locally compact group is relatively compact if
and only if it is \emph{locally bounded}, i.e., there is a neighborhood
$U_{e_G}$ whose image is contained in some element of the filter $\mathfrak H$
of neighborhoods of $e_H$ in~$H$ having compact closure.

Let us also recall the notion of discontinuity group of a homomorphism $\pi$ of
a topological group $G$ into a topological group $H$, which was introduced
in~[1] and used in~[2, 3].

\begin{definition}\label{d1}
Let $G$ and $H$ be topological groups, and let $\mathfrak U=\mathfrak U_G$ be
the filter of neighborhoods of $e_G$ in~$G$. For every (not necessarily
continuous) locally relatively compact homomorphism $\pi$ of~$G$ into~$H$, the
set $\operatorname{DG}(\pi)=\bigcap_{U\in\mathfrak U}{\overline{\pi(U)}} $ is
called the discontinuity group of~$\pi$.
\end{definition}

Here and below, the bar stands for the closure in the corresponding topology
(therefore, here the closure is taken with respect to the topology of~$H$).
(See~[2, Definition~1.1.1].)

The discontinuity group of a homomorphism has some important properties. Let us
list some of them.

Let $\pi$ be a (not necessarily continuous) locally relatively compact
homomorphism of a connected Lie group~$G$ into a topological group~$H$. Then
the set $\operatorname{DG}(\pi)$ is a compact subgroup of the topological
group~$H$ and a compact normal subgroup of the closed
subgroup~$\overline{\pi(G)}$ of~$H$.

Moreover, the filter basis $\{{\overline{\pi(U)}}\mid U\in\mathfrak U\}$
converges to~$\operatorname{DG}(\pi)$, and the homomorphism~$\pi$~is continuous
if and only if $\operatorname{DG}(\pi)=\{e_H\}$. (See~[2, Theorem~1.1.2].)

Under the same conditions,~$\operatorname{DG}(\pi)$ is a compact connected
subgroup of~$H$. (See~[2, Lemma~1.1.6].)

Let~$G$~be a connected locally compact group, let $N$~be a closed normal
subgroup of~$G$, and let~$\pi$ be a locally relatively compact homomorphism
of~$G$ into a topological group~$H$ \textup(for instance, a locally bounded
homomorphism into a locally compact group\/\textup). Let~$M$ be the
discontinuity group of the restriction~$\operatorname{DG}(\pi|_N)$. Then~$M$ is
a closed normal subgroup of the compact discontinuity
group~$\operatorname{DG}(\pi)$, and the corresponding quotient
group~$\operatorname{DG}(\pi)/M$ is isomorphic to the discontinuity group
~$\operatorname{DG}(\psi)$ of the homomorphism~$\psi$ of~$G$ obtained as the
composition of the homomorphism~$\pi$ and the canonical
homomorphism~$\overline{\pi(G)}\to\overline{\pi(G)}/M$. (See~[2, Lemma~1.1.7].)

In this note we prove that the discontinuity group of every locally bounded
homomorphism of a Lie group into a Lie group is not only compact and connected
but also commutative.

\section{Preliminaries}

The following refinement of Lemma~1.1.9 of~[2] obviously holds.

\begin{theorem}\label{th1} Let $G$ be a topological group, let~$G'$ be the commutator
subgroup of~$G$ \textup(in the algebraic sense and in the topology induced by
that of~$G$\/\textup), and let~$\pi$ be a locally relatively compact
homomorphism of~$G$ into some topological group~$H$. Then the commutator
subgroup of the discontinuity group of the homomorphism~$\pi$ is a subgroup of
the discontinuity group of the restriction~$\pi|_{G'}$ of~$\pi$ to the
commutator subgroup~$G'$ of~$G$\textup:
$$
\operatorname{DG}(\pi)'\subset\operatorname{DG}(\pi|_{G'});
$$
moreover, $\operatorname{DG}(\pi)'$ is a normal subgroup
of~$\operatorname{DG}(\pi|_{G'})$.
\end{theorem}

\textbf{Proof.} The group $\operatorname{DG}(\pi)'$ is a normal subgroup
of~$\operatorname{DG}(\pi|_{G'})$ indeed since, by the invariance of $G'$
in~$G$, $\operatorname{DG}(\pi)'$ is obviously invariant with respect to inner
automorphisms of~$G$, and hence normal in~$\operatorname{DG}(\pi|_{G'})$.
\smallskip

The proof of Lemma~1.1.9 in~[2] is a model for the proof of the following main
theorem.

\section{Main result}

\begin{theorem}\label{t2}
 The discontinuity group of every locally bounded homomorphism of
a Lie group into a Lie group is commutative.
\end{theorem}

\textbf{Proof.} Let $G$ and $H$ be Lie groups, and let $\pi$ be a locally
bounded homomorphism of $G$ into $H$. Let $G^9$ be the connected component of
$G$. The discontinuity group of $\pi$ coincides obviously with the
discontinuity group of the restriction $\pi|_{G^0}$. Thus, we may assume that
$G$ is connected. Let $DG(\pi)\subset H$ be the discontinuity group of~$\pi$.

Consider the adjoint representation $\operatorname{Ad}H$ of $H$. The
composition $\operatorname{Ad}H\circ\pi$ is a locally bounded linear
representation of~$G$. By Corollary~1.3.3 of~[2], the discontinuity
group~$DG(\operatorname{Ad}H\circ\pi)$ is a commutative compact connected
subgroup of~$H$.

By the very definition $\operatorname{DG}(\pi)=\bigcap_{U\in\mathfrak
U}{\overline{\pi(U)}}$ of the discontinuity group of a homomorphism and by the
continuity of $\operatorname{Ad}H$, we see that, since the kernel of
$\operatorname{Ad}H$ coincides with the center $Z_H$ of the group~$H$, it
follows that
$\operatorname{DG}(\operatorname{Ad}H\circ\pi)=\bigcap_{U\in\mathfrak
U}{\overline{(\operatorname{Ad}H\circ\pi)(U)}}$ coincides with
$\operatorname{Ad}H(DG(\pi))$. Therefore,
$\operatorname{DG}(\operatorname{Ad}H\circ\pi)$ is the image
of~$\operatorname{DG}(\pi)$ under the factorization of $H$ by the center~$Z_H$.

However, both $Z_H$ and~$\operatorname{DG}(\operatorname{Ad}H\circ\pi)$ are
commutative groups. Therefore,~$\operatorname{DG}(\pi)$ is solvable. Let $Z_H^0$
be the connected component of the identity element $e$ of~$Z_H$ and let $D$ be
the complementary discrete subgroup of~$Z_H$ to~$Z_H^0$ such that $D\cap
Z_H^0=\{e\}$. Then the connected component of~$\operatorname{DG}(\pi)$ is a
connected solvable Lie subgroup of the compact Lie
group~$\operatorname{DG}(\pi)$. Hence, this component is commutative. The whole
group~$\operatorname{DG}(\pi)$ is generated by this commutative connected
component of the identity element and the central subgroup~$D$.
Therefore,~$\operatorname{DG}(\pi)$ is commutative, as was to be proved.
\smallskip

\section{Discussion}

As is known, Theorem~2 cannot be extended even to connected compact groups,
since they can have discontinuous non-one-dimensional finite-dimensional
representations~[1]--[3].

\section*{Funding}

The research was supported by SRISA RAS according to the project FNEF-2022-0007
(Reg. no.~1021060909180-7-1.2.1).

\end{document}